\documentclass[10pt]{amsart}

\usepackage{sudokuLT}
\usepackage{shidokuLT}
\usepackage{amssymb,amsthm,amsmath}
\usepackage[curve]{xypic}

\newtheorem{theorem}{Theorem}[]

\theoremstyle{remark}

\numberwithin{equation}{section}
\theoremstyle{theorem}

\newtheorem*{fixinglemmas}{Fixing Lemmas}

\usepackage[english]{babel}
\usepackage[pdftex]{graphicx}
\usepackage{color}

\usepackage{amssymb,amsthm,amsmath}

\usepackage{epsfig}
\usepackage{verbatim}

\renewcommand*\shidokuformat[1]{\small $#1$} 
\renewcommand*\sudokuformat[1]{\small $#1$} 
\begin{document}

\title{Minimal complete Shidoku symmetry groups}

\author{Elizabeth Arnold, Rebecca Field, Stephen Lucas, and Laura Taalman}
\address{Department of Mathematics and Statistics, MSC 1911, James Madison University, Harrisonburg, VA 22807}
\email{
{\texttt{[arnoldea,fieldre,lucassk,taalmala]@math.jmu.edu}}}
\urladdr{\texttt{http://www.math.jmu.edu/\~[arnoldea,fieldre,lucassk,taalmala]/}}

\date{October 13, 2011}

\maketitle

\begin{abstract}

Calculations of the number of equivalence classes of Sudoku boards has to this point been done only with the aid of a computer, in part because of the unnecessarily large symmetry group used to form the classes.  In particular, the relationship between relabeling symmetries and positional symmetries such as row/column swaps is complicated.  In this paper we focus first on the smaller Shidoku case and show first by computation and then by using connectivity properties of simple graphs that the usual symmetry group can in fact be reduced to various minimal subgroups that induce the same action.  This is the first step in finding a similar reduction in the larger Sudoku case and for other variants of Sudoku.

\end{abstract}


\section{Introduction: Sudoku and Shidoku}

A growing body of mathematical research has demonstrated that one of the most popular logic puzzles in the world, Sudoku, has a rich underlying structure. 
In Sudoku, one places the digits from one to nine in a nine by nine grid with the constraint that no number is repeated in any row, column, or smaller three by three block. We call these groups of nine elements \textit{regions}, and the completed grid a \textit{Sudoku board}. A subset of a Sudoku board that uniquely determines the rest of the board is a Sudoku puzzle, as illustrated in Figure 1, taken from \cite{nofrills}.  Note that this puzzle contains only eighteen starting clues, which is the conjectured minimum number of clues for a 180-degree rotationally symmetric Sudoku puzzle.

\begin{figure}[h,t]
\begin{center}
\begin{sudoku-block} 
| | | | | | | |5|7|. 
| | | |8| | | | |2|. 
| | | |3|9| | | | |. 
| | | | | |1|6| | |. 
| | |8| | | |9| | |. 
| | |2|7| | | | | |. 
| | | | |5|2| | | |. 
|1| | | | |7| | | |. 
|8|6| | | | | | | |. 
\end{sudoku-block} 
\hspace{.5pc}
\begin{sudoku-block} 
|4|8|9|2|1|6|3|5|7|.
|5|3|6|8|7|4|1|9|2|.
|2|7|1|3|9|5|8|6|4|.
|7|5|3|9|4|1|6|2|8|.
|6|4|8|5|2|3|9|7|1|.
|9|1|2|7|6|8|5|4|3|.
|3|9|4|1|5|2|7|8|6|.
|1|2|5|6|8|7|4|3|9|.
|8|6|7|4|3|9|2|1|5|.
\end{sudoku-block} 
\end{center}
\vspace{-.25pc}

\caption{An 18-clue Sudoku puzzle and its unique solution board.}
\end{figure}

The current body of work on Sudoku ranges from popular interest monographs such as Wilson \cite{Wilson} and appeals to the mathematical interest in Sudoku and its variants for an undergraduate math audience as in Taalman \cite{Taalman} to a wide range of more theoretical papers.  For example Bailey, Cameron and Connelly \cite{gerechte} explore mutually orthogonal collections of Sudoku boards, Arnold \textit{et al.}  \cite{Arnold} and Sato \textit{et al.} \cite{Sato} apply Gr\"obner basis techniques to Sudoku, and Newton and DeSalvo \cite{Newton} investigate the Shannon entropy of collections of Sudoku matrices.  Most relevant to our current research is a sophisticated piece of computation by Felgenhauer and Jarvis \cite{Felgenhauer} showing that there are 6,670,903,752,021,072,936,960 different Sudoku boards and a related paper by Russell and Jarvis \cite{Russell} that showed that these boards split into 5,472,730,538 equivalence classes under the action of the standard Sudoku symmetry group. 

In this paper, we consider various symmetry groups for the more accessible variant of Sudoku known as Shidoku. A {\em Shidoku board} is a $4 \times 4$ Latin square whose regions (rows, columns, and designated $2 \times 2$ blocks) each contain the integers one to four exactly once.  In this smaller universe it is not too difficult to show that there are 288 different Shidoku boards \cite{Taalman}.   There are a small number of generating symmetries for Shidoku and while the full Shidoku symmetry group is well known to be relatively large, we establish several different subgroups of minimal size which generate the same equivalence relation among Shidoku boards as the full group of symmetries.  We follow this with a discussion of Burnside's lemma (as used in \cite{Russell}) where we demonstrate one practical application of such minimal symmetry groups.  
Finally, we show how all minimal symmetry groups satisfying various conditions may be attained by considering equivalence classes of group actions on graphs.  The Russell and Jarvis computation \cite{Russell} involved computer methods in part because of the large size of the Sudoku symmetry group.  Our work is a first step in simplifying such a calculation by reducing to a smaller, but still complete Sudoku symmetry group. 

\section{The Full Shidoku Symmetry Group}

A {\em Shidoku symmetry} is a map from the set $\mathcal{S}$ of Shidoku boards to itself.  We will call elements of $S_4$ that permute the values $\{1,2,3,4\}$ on the board {\em relabeling symmetries}, and elements of $S_{16}$ that permute the cells of the board while preserving the Shidoku conditions on every possible board {\em position symmetries}.  Note that every element of $S_4$ defines a valid relabeling symmetry, but not every element of $S_{16}$ defines a valid position symmetry.  For example, swapping the first and second rows of any Shidoku board preserves the Shidoku conditions, but swapping the second and third rows does not always preserve the Shidoku conditions. 

Every Shidoku symmetry is a combination of the set $S_4$ of relabeling symmetries and the set $H_4$ of position symmetries described above.  Therefore the {\em full Shidoku symmetry group} is the direct product $G_4 = H_4 \times S_4$.  The action of this group partitions the set $\mathcal{S}$ of Shidoku boards into two orbits, and we call two Shidoku boards are {\em equivalent} if they are in the same orbit under the action of the full Shidoku symmetry group $G_4$.  The 288 Shidoku boards split into two such equivalence classes, one with 96 boards which we refer to as {\em Type 1} boards, one of which is on the left in Figure 2, and one with 192 boards, called {\em Type 2} boards, one of which is given on the right in Figure 2 (\cite{Arnold}, \cite{Taalman}).

\begin{figure}[h]
\begin{center}
\begin{shidoku-block}
|1|2|3|4|.
|3|4|1|2|.
|2|1|4|3|.
|4|3|2|1|.
\end{shidoku-block}
\qquad
\begin{shidoku-block}
|1|2|3|4|.
|3|4|1|2|.
|2|3|4|1|.
|4|1|2|3|.
\end{shidoku-block}
\end{center}
\caption{Type 1 and Type 2 Shidoku representatives.}
\end{figure}
\vspace{-.25pc}

The position symmetry group $H_4$ is generated by the following symmetries, where a {\em band} is the combination of either the first and second or the third and fourth rows, and a {\em pillar} is an analogous combination of columns (see \cite{Felgenhauer}):\\*[-.5pc]
\begin{itemize}
\item swapping rows/columns within bands/pillars \\*[-.3pc]

\item swapping bands/pillars \\*[-.3pc]

\item rotation $r$ of the board by a quarter turn clockwise \\*[-.3pc]

\item transpose $t$ of the board\\*[-.4pc]
\end{itemize}

Fortunately, we can obtain a much smaller set of generators for $H_4$.  In fact, we need only rotation, transpose, and one row swap:  let $s$ to be the swap of the third and fourth rows of a Shidoku board, then $H_4$ is generated by $r$, $s$, and $t$.  (For example, to swap the first and second columns we can apply $r^{-1}$, then $s$, then $r$.)   There are relations among these three generators; the easiest to see being the orders of the generators and a few other relations, such as 
$trtr=id$.  The full presentation of the position symmetry group is
$$H_4 = \langle r, s, t \mid r^4, s^2, t^2, trtr, rsr^2tsr^3t, tstststs, srsr^3srsr^3\rangle.$$

\noindent
We will use the more compact notation $H_4 = \langle r, s, t \rangle$.  This is a non-abelian group of order 128 (see \cite{GAP,lorch}), so the full Shidoku symmetry group $G_4 = H_4 \times S_4$ has order $128(4!)=3072$.

Note that we have now defined a group of relatively large order, $|G_4|=3072$, acting on a set of much smaller order, $|\mathcal{S}|=288$, to make two orbits, the largest of which is size $192$.  We will say that a group of Shidoku symmetries is {\em complete} if its action partitions the set of Shidoku boards into the same Type 1 and Type 2 orbits.  The natural question to ask is whether we can find a complete group of Shidoku symmetries with order smaller than $3072$?  The minimum possible order that such a group could have is the size $192$ of the largest orbit,   so can we find a complete Shidoku symmetry group with this minimum order?

We will start with a heuristic investigation using {\em MATLAB} \cite{Matlab}, {\em GAP} \cite{GAP}, and the graph visualization program {\em yEd} \cite{yed}. The action of the full group $G_4 = H_4 \times S_4$ on the set $\mathcal{S}$ of Shidoku boards can be represented with the graph in Figure 3, where each vertex is one of the $288$ Shidoku boards and each edge corresponds to one of the generators $r$, $s$, $t$, $(1 2)$, $(2 3)$, $(3 4)$, and $(1 4)$.  (We use more generators of $S_4$ than necessary for  symmetry.)
\vspace{1pc}

\begin{figure}[h]
\centerline{\includegraphics[height=4in, angle=90]{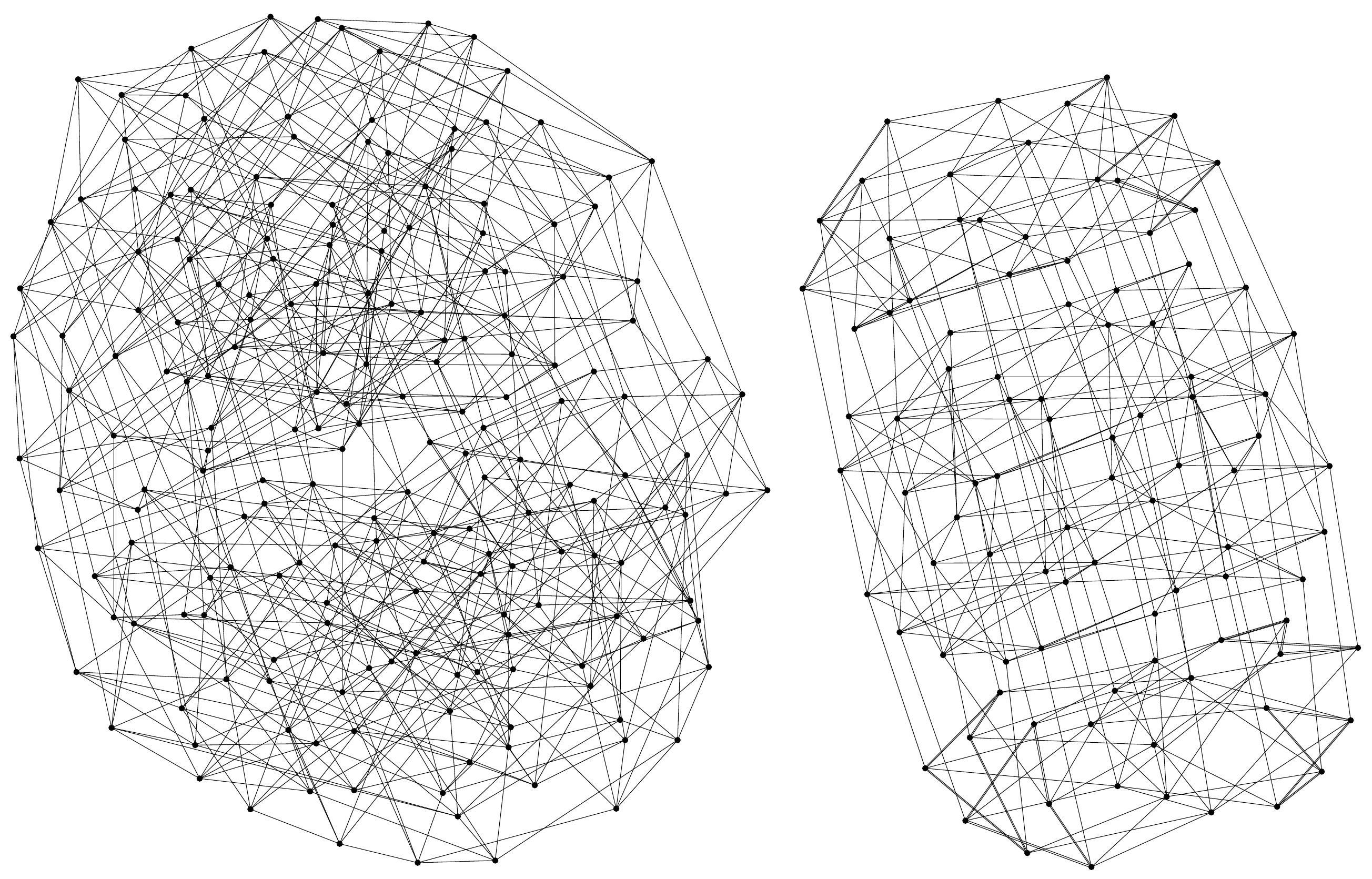}}
\caption{The action of $G_4=H_4 \times S_4$ on $\mathcal{S}$.}
\end{figure}

This graph clearly shows the smaller Type 1 orbit (on the right) and the larger Type 2 orbit (on the left).  Given the high degree of connectivity in this graph, it is reasonable to hope that we can remove the edges for one or more generators and still retain a two-component graph.

\section{Minimal Complete Shidoku Symmetry Groups}

Via a sequence of examples we now show that it is possible to construct minimal complete Shidoku symmetry groups.  We will return to this in Section \ref{S:vizS4} from a more algebraic perspective that is better suited for future generalization to the $9 \times 9$ Sudoku case.

\begin{theorem}
There exist complete Shidoku symmetry groups of minimal order $192$.
\label{T:minimal}
\end{theorem}

Notice that since the minimum possible order of any complete Shidoku symmetry group is $192$, any complete Shidoku symmetry group will need to include a mix of position and relabeling symmetries as the groups $S_4$ and $H_4$ are each too small on their own.
\vspace{1pc}

\noindent
{\em Example 1:}
$\langle r, t \rangle \times S_4$

\noindent
Removing the row swap $s$ leaves the 
remaining position symmetry group denoted $\langle r, t \rangle$.  This is the order 8 subgroup $\langle r, t \mid r^4, t^2, trtr \rangle \subseteq H_4$.  Therefore $\langle r, t \rangle \times S_4$ has order 192, precisely the minimum order we wish to obtain.  Unfortunately, as seen in Figure 4, the action of this group has five orbits and therefore this subgroup is not a complete
Shidoku symmetry group.
\vspace{1pc}

\begin{figure}[h,t]
\centerline{\includegraphics[height=2.8in]{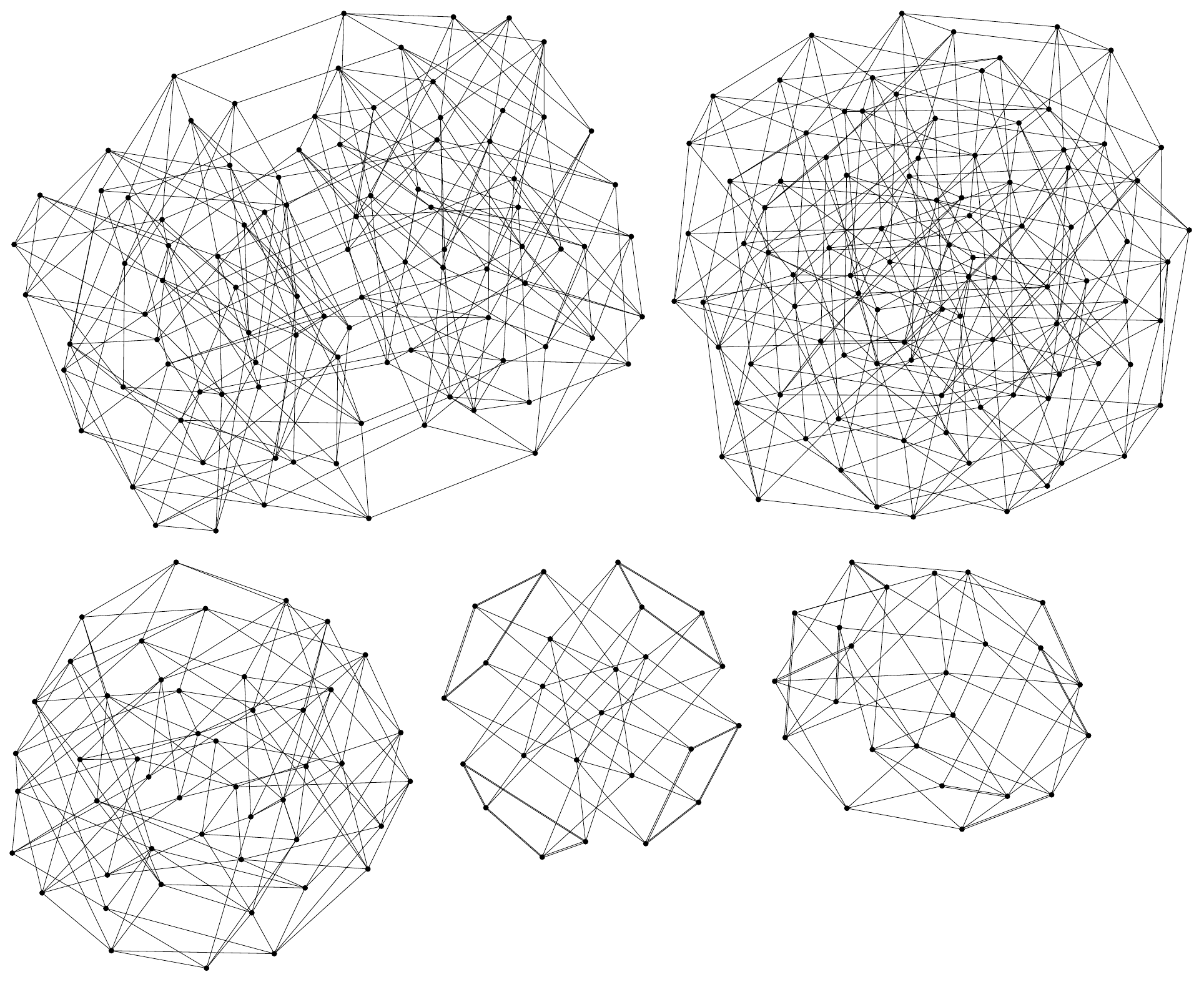}}
\caption{The action of $\langle r, t \rangle \times S_4$ on $\mathcal{S}$.}
\end{figure}

\noindent
{\em Example 2:}
$\langle r, s \rangle \times S_4$ and $\langle r, s \rangle \times \langle (1 2 3) \rangle$

\noindent
Removing the transpose symmetry $t$ leaves the position symmetry subgroup $\langle r, s \rangle \subseteq H_4$.  Since this subgroup has order $64$, 
we can obtain a complete Shidoku symmetry group $\langle r, s \rangle \times S_4$ of order $64(4!)=1536$.  
A group fo minimal size can be obtained by extending $\langle r,s\rangle$ by 
an order three subgroup of $S_4$, 
for example $\langle (1 2 3) \rangle \subseteq S_4$, and this does combine with $\langle r, s \rangle$ to form a minimal complete Shidoku symmetry group, as shown in Figure 5.  Note the somewhat surprising corollary that any valid Shidoku board can be transformed into any other equivalent board using nothing more than rotation, a single row swap, and a single cyclic permutation of order three.  
\vspace{1pc}

\begin{figure}[h,t]
\centerline{\includegraphics[height=2in]{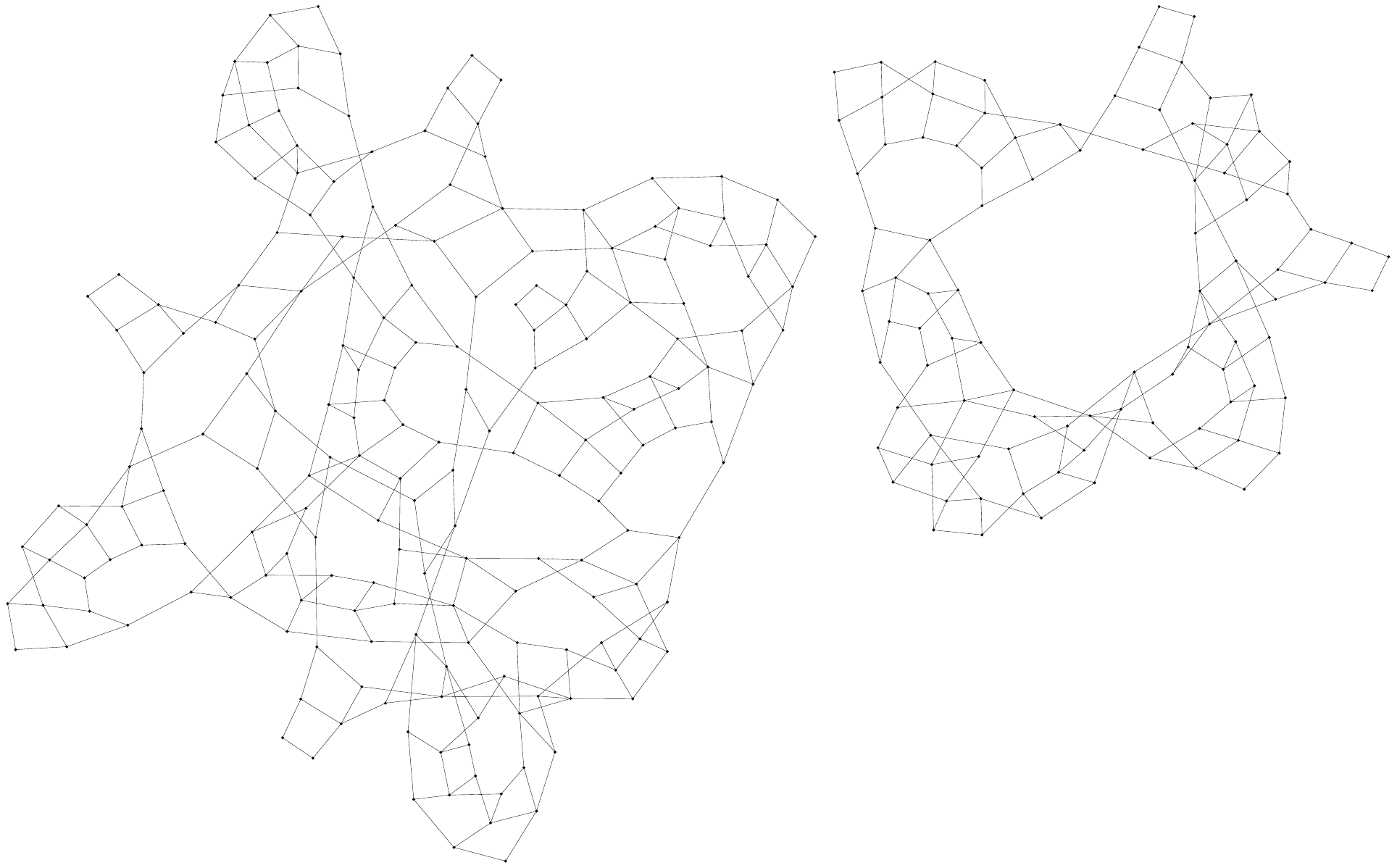}}
\caption{The action of $\langle r, s \rangle \times \langle (1 2 3) \rangle$ on $\mathcal{S}$.}
\end{figure}

\noindent
{\em Example 3:}
$\langle s, t \rangle \times S_4$

\noindent
In Examples 2 and 3, we saw that removing the swap symmetry $s$ from $H_4$ results in a Shidoku symmetry group that is minimal but not complete, and removing the transpose symmetry $t$ results in a Shidoku symmetry group that is complete but not minimal (unless we also pass to a subgroup of $S_4$).  Removing the rotation symmetry $r$ gives us the best of both worlds.  The position symmetry subgroup $\langle s, t \rangle$ has order 8 so  $\langle s, t \rangle \times S_4$ has order 192 and therefore is minimal.  As shown in Figure 6, this group induces two orbits on $\mathcal{S}$ and therefore is complete.
\vspace{1pc}

\begin{figure}[h,t]
\centerline{\includegraphics[height=2in]{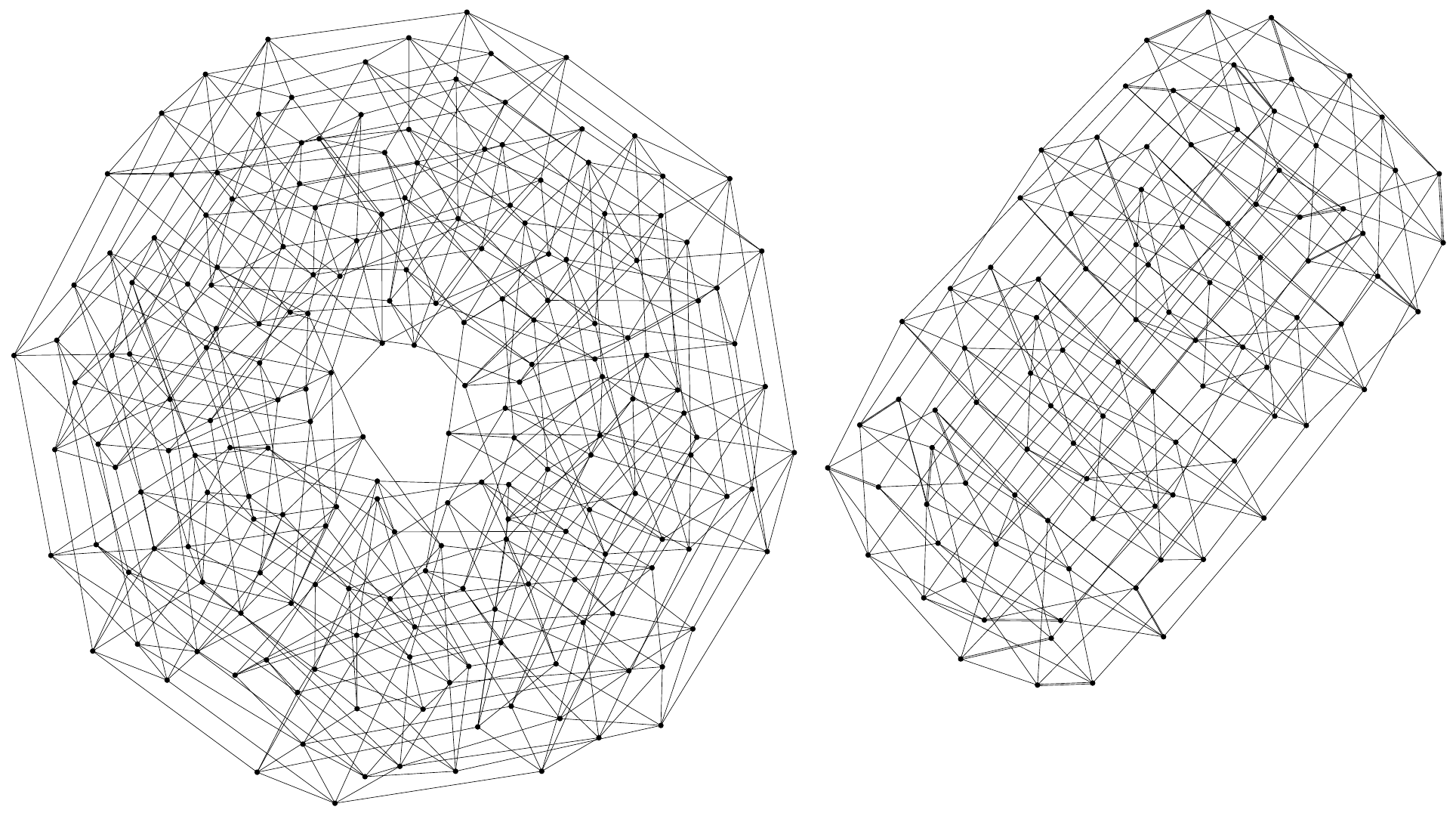}}
\caption{The action of $\langle s, t \rangle \times S_4$ on $\mathcal{S}$.}
\end{figure}

\noindent
{\em Example 4:}
$H_4 \times \langle (1 2 3) \rangle$ and $\langle r^2, s, t \rangle \times \langle (1 2 3) \rangle$

\noindent
Finally, we try passing to a subgroup of $S_4$ while keeping all of the position symmetries.  The full position symmetry group $H_4$ has order 128, and 384 is the smallest multiple of 128 that is divisible by 192.  Therefore, a complete Shidoku symmetry group containing all of $H_4$ must contain a subgroup of $S_4$ whose order is a multiple of $3$.  An obvious candidate is $\langle (1 2 3) \rangle \in S_4$.  Combining these relabelings with all of $H_4$ does result in a graph with the two desired orbits, as shown in Figure 7.

\begin{figure}[h]
\centerline{\includegraphics[height=2in]{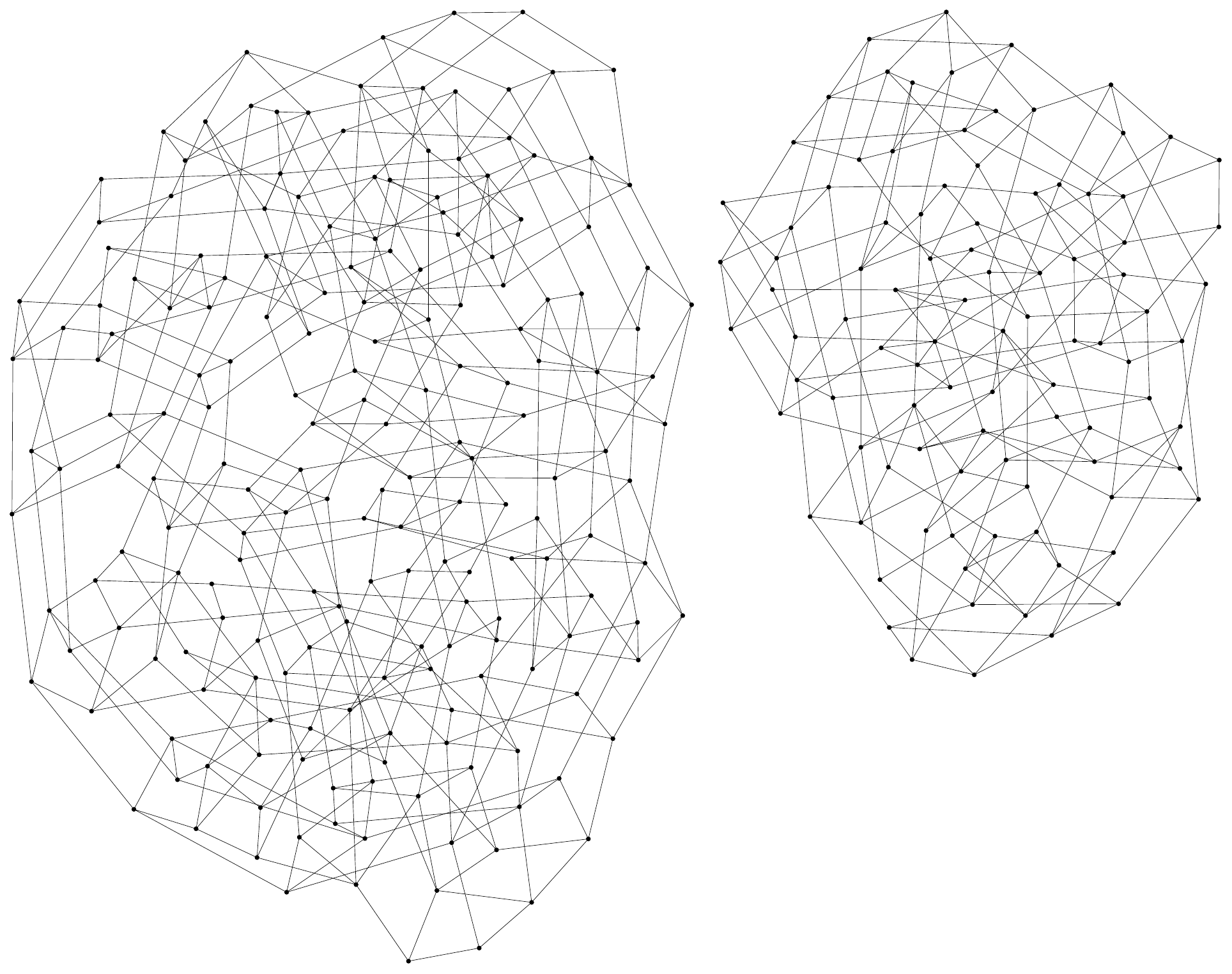}}
\caption{Figure 7: The action of $H_4 \times \langle (1 2 3) \rangle$ on $\mathcal{S}$.}
\end{figure}

Although $H_4 \times \langle (1 2 3) \rangle$ is a complete Shidoku symmetry group, it has order 384 and therefore is not minimal.  However, we saw in Example 3 that the rotation symmetry $r$ is, in some sense, redundant.  Namely, given any board $A \in \mathcal{S}$, the board $r(A)$ can be obtained by the composition of $s, t$ and a relabeling symmetry that depends on the board $A$. Replacing $r$ (rotation by 90 degrees) with $r^2$ (rotation by 180 degrees) results in the complete minimal symmetry group $\langle r^2, s, t \rangle \times \langle (1 2 3) \rangle$.  


\section{Using Burnside's Lemma to Count Orbits}

The number of equivalence classes of Sudoku boards can be calculated by applying Burnside's Lemma, which states that if a finite group acts on a set, then the number of orbits is the average size of the fixed point sets for the elements of the group.  For Shidoku, this means that the number of orbits under the action of the group of Sudoku symmetries is the average of the number of Sudoku boards fixed by each symmetry group element. 

Russell and Jarvis applied Burnside's Lemma to the full group of all possible Sudoku symmetries in \cite{Russell} using a computer.  Their brute-force computation found 5,472,730,538 equivalence classes of Sudoku boards via $275$ conjugacy classes in a Sudoku position symmetry group of size $3359232(9!)$.  If a smaller but complete group of Sudoku symmetries were known, then this calculation could be greatly simplified.  For Shidoku, the reduction of the symmetry group from $H_4 \times S_4$ to the minimal $\langle s, t \rangle \times S_4$ will allow us to easily perform this calculation by hand.  
As seen in Figure 6, there 
are two orbits of Shidoku boards under the action of this subgroup; our goal here is to illustrate how tractable the problem becomes when we work with a minimal Shidoku symmetry group, and thus motivate the search for such groups.

We will say that a Shidoku board $B$ is {\em invariant} under a position symmetry $x$ if we can undo the action of $x$ with a relabeling symmetry; that is, if there exists a relabeling symmetry $\sigma$ so that $\sigma(x(B)) = B$.  If a board $B$ is invariant under a symmetry $x$ then $x(B)$ differs from $B$ by a relabeling. 

More precisely, a Shidoku board $B$ is an ordered list of values $(b_1,b_2,b_3,\ldots,b_{16})$, reading left to right and then top to bottom.  A Shidoku symmetry $x$ is then an element of $S_{16}$ that produces another valid Shidoku board.  For each $i$, a symmetry $x$ takes the value $b_i$ that is in the $i$th cell of $B$ and moves it to the $x(i)$th cell of $x(B)$.
\vspace{1pc}

\noindent
{\em Example 5: Invariance under transpose}

\noindent 
Let $t$ be the position symmetry taking the transpose of a Shidoku board.  In permutation notation, $t = (2\,\,\, 5)(3\,\,\,9)(4\,\,\,13)(7\,\,\,10)(8\,\,\,14)(12\,\,\,15)$.  The action of $t$ is pictured on the right in Figure 8.  Notice that, for example, the value $b_8$ in the 8th cell of $B$ moves to the $t(8)=14$th cell of $t(B)$.

\begin{figure}[h]
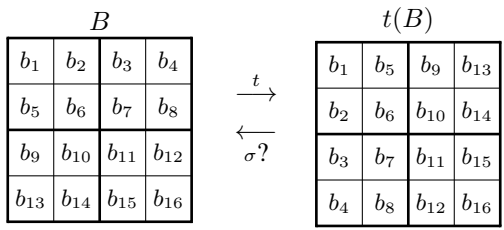

\begin{center}
\parbox{1.3in}{
\centerline{$B$}
\begin{shidoku}
|{\mbox{$b_1$}}|{\mbox{$b_2$}}|{\mbox{$b_3$}}|{\mbox{$b_4$}}|.
|{\mbox{$b_5$}}|{\mbox{$b_6$}}|{\mbox{$b_7$}}|{\mbox{$b_8$}}|.
|{\mbox{$b_9$}}|{\mbox{$b_{10}$}}|{\mbox{$b_{11}$}}|{\mbox{$b_{12}$}}|.
|{\mbox{$b_{13}$}}|{\mbox{$b_{14}$}}|{\mbox{$b_{15}$}}|{\mbox{$b_{16}$}}|.
\end{shidoku}}
\hspace{-2.15in}
\parbox{.4in}{
$$\overset{t}{\longrightarrow}$$
$$\underset{\sigma\mbox{?}}{\longleftarrow}$$}
\parbox{1.3in}{
\centerline{$t(B)$}
\begin{shidoku}
|{\mbox{$b_1$}}|{\mbox{$b_5$}}|{\mbox{$b_9$}}|{\mbox{$b_{13}$}}|.
|{\mbox{$b_2$}}|{\mbox{$b_6$}}|{\mbox{$b_{10}$}}|{\mbox{$b_{14}$}}|.
|{\mbox{$b_3$}}|{\mbox{$b_7$}}|{\mbox{$b_{11}$}}|{\mbox{$b_{15}$}}|.
|{\mbox{$b_4$}}|{\mbox{$b_8$}}|{\mbox{$b_{12}$}}|{\mbox{$b_{16}$}}|.
\end{shidoku}}
\end{center}
\caption{If $\sigma \in S_4$ exists then board $B$ is invariant under $t$.}
\end{figure}

\noindent
A board $B$ is invariant under $t$ if there exists some relabeling $\sigma$ that returns $t(B)$ to $B$.  In general, a Shidoku board $B$ is invariant under a symmetry $x$ if there is some relabeling permutation $\sigma$ so that for each $i$ we have $b_{x(i)} = \sigma(b_i)$.  Notice that $\sigma$ will depend on values of the $b_i$; for example, in this case, since $t(8)=14$, we must have $b_{t(8)}=\sigma(b_8)$.   
Figure 9 illustrates that for this particular board $B$, the relabeling that returns $t(B)$ to the original board $B$ is $\sigma = (2\,\,\,3)$.

\begin{figure}[h]
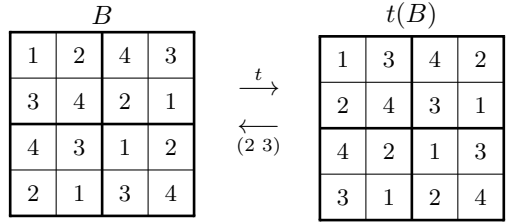

\begin{center}
\parbox{1.3in}{
\centerline{$B$}
\begin{shidoku}
|1|2|4|3|.
|3|4|2|1|.
|4|3|1|2|.
|2|1|3|4|.
\end{shidoku}}
\hspace{-2.15in}
\parbox{.4in}{
$$\overset{t}{\longrightarrow}$$
$$\underset{(2\,\,3)}{\longleftarrow}$$}
\parbox{1.3in}{
\centerline{$t(B)$}
\begin{shidoku}
|1|3|4|2|.
|2|4|3|1|.
|4|2|1|3|.
|3|1|2|4|.
\end{shidoku}}
\end{center}
\caption{Board $B$ is invariant under the action of $t$, up to relabeling.}
\end{figure}

Let us now revisit the minimal Shidoku symmetry group $\langle s, t \rangle \times S_4$.  The subgroup\\ $\langle s, t \rangle = \{id, s, t, st, ts, sts, tst, stst\} \subseteq H$ is the Coxeter group with presentation $\langle s, t \mid s^2, t^2, (st)^4 \rangle$, where the last relation is equivalent to the relation $s(tst)=(tst)s$ that follows from the fact that the row swap $s$ commutes with the column swap $tst$.
This subgroup partitions into five conjugacy classes, with representatives $id$, $s$, $t$, $st$, and $stst$.  Within each class, each position symmetry has the same number of invariant boards.  The number of invariant Shidoku boards for each element of each conjugacy class is given in Table 1.
\vspace{.5pc}

\begin{table}[h]
\begin{center}
\begin{tabular}{llr}
\multicolumn{1}{l}{Class}
  & \multicolumn{1}{l}{Representative}
  & \multicolumn{1}{l}{Invariant} \\
\hline
\hline
$C_{id} = \{id\}$ & $()$ & $12 \cdot 4!$ \\
\hline
$C_s = \{s, tst\}$ & $(9\,\,13)(10\,\,\,14)(11\,\,15)(12\,\,\,16)$ & $0 \cdot 4!$ \\
\hline
$C_t = \{t, sts\}$ & $(2\,\,\,5)(3\,\,\,9)(4\,\,13)(7\,\,10)(8\,\,14)(12\,\,\,15)$ & $2 \cdot 4!$ \\
\hline
$C_{st} = \{st, ts\}$ & $(2\,\,\,5)(3\,\,\,9\,\,\,4\,\,\,13)(7\,\,\,19\,\,\,8\,\,\,4)(11\,\,12\,\,\,16\,\,\,15)$ & $0 \cdot 4!$ \\
\hline
$C_{stst} = \{stst\}$ & $(3\,\,\,4)(7\,\,\,8)(9\,\,\,13)(10\,\,\,14)(11\,\,16)(12\,\,\,15)$ & $0 \cdot 4!$ \\
\hline
\end{tabular}
\end{center}
\caption{Number of invariant boards for each conjugacy class of $\langle s, t \rangle$.}
\end{table}

Note that the identity element fixes all $12 \cdot 4!$ Shidoku boards.  On the other hand, the element $s$ fixes no Shidoku boards up to permutation; that is, there is no board $B$ for which a relabeling symmetry can undo the action of the row swap $s$.   Now consider the transpose symmetry $t$ in the third row of Table 1.  In Example 5 we saw that a Shidoku board can only be invariant under this symmetry if the returning relabeling fixes the values on the main diagonal and swaps the remaining two values in the upper-left block.  It is easy to show that for any given labeling of elements in the upper-left block there are only two Shidoku boards that are invariant under transposition; thus there are $2 \cdot 4!$ Shidoku boards that are invariant under the action of $t$.

In general, in order to construct a Burnside's table for a given Shidoku or Sudoku symmetry group, we need to consider various relationships between position and relabeling symmetries such those in the following ``Fixing Lemmas.''

\begin{fixinglemmas}
Suppose $B$ is a Shidoku board that is invariant under the action of a Sudoku symmetry $x$ via a relabeling permutation $\sigma$.  Then we have the following relationships:

\begin{enumerate}

\item Value-Fixing:$\;$  If $x$ fixes a cell whose value is $n$, then $\sigma$ must fix that $n$.  In other words, if $x(i)=i$ and $b_i=n$ then $\sigma(n)=n$.

\item Block-Fixing:$\;$  If $x$ fixes any entire row, column, or block of $B$, then $\sigma$ must be the identity relabeling.

\item Fixed Points:$\;$  If $\sigma$ fixes $n$, then $x$ must take every cell of $B$ whose value is $n$ to another cell whose value is $n$. In other words, if $\sigma(n)=n$ and $b_i=n$, then $b_{x(i)}=n$.

\end{enumerate}
\label{L:fixing}
\end{fixinglemmas}

\begin{proof}
Let $B$ be a Shidoku or Sudoku board, $x$ a Shidoku or Sudoku symmetry, and $\sigma$ a relabeling symmetry such that $\sigma(x(B))=B$.  To prove (i), suppose $x(i)=i$ and $b_i=n$ for some index $i$ and value $n$.  Combining this information with the definition of invariance we have
$$\sigma(n) = \sigma(b_i) = b_{x(i)} = b_i = n.$$

Since every row, column, or block region contains all of the numbers 1--4 we immediately obtain (ii) as a corollary to (i).

Finally, to prove (iii), suppose $\sigma(n)=n$ and $b_i=n$.  Combining this information with the definition of invariance we have
$$b_{x(i)} = \sigma(b_i) = \sigma(n) = n.$$

\vspace{-1.5\baselineskip}
\end{proof}

There are many other ``fixing lemmas'' that we can use to develop techniques for counting invariant boards for other symmetry groups.  Applying these techniques to the full Shidoku symmetry group $H_4 \times S_4$ we can produce a table of invariant board counts similar to Table 1, although with twenty rows of conjugacy classes instead of just five.

\vspace{1pc}
Both $H_4 \times S_4$ and $\langle s, t \rangle \times S_4$ are complete symmetry groups, and therefore have two orbits in $\mathcal{S}$.  Thus by Burnside's Lemma, the average over all group elements of the number of invariant boards must be two.  Applying Burnside's Lemma to the invariance data in Table 1 we see that for the minimal Shidoku symmetry group $\langle s, t \rangle \times S_4$ the number of orbits is indeed
$$\frac{1(12 \cdot 4!) + 2(0) + 2(2 \cdot 4!) + 2(0) + 1(0)}{8\cdot 4!} = 2.$$

For $4 \times 4$ Shidoku we can do these invariance calculations by hand.  For $9 \times 9$ Sudoku there are $275$ conjugacy classes of the full position Sudoku symmetry group of order 3,359,232 (see \cite{Russell}).  In this case finding the conjugacy classes, the invariance calculations, and even the final Burnside's Lemma computation must be done by computer.  Reducing the size of the Sudoku symmetry group and extending the Fixing Lemma techniques could pave the way to a more straightforward non-computer calculation of the number of equivalence classes of $9 \times 9$ Sudoku boards.

\vspace{1pc}

All of the above information on invariant boards and orbits is contained in the graph for $\langle s, t \rangle \times S_4$ acting on $\mathcal{S}$ from Figure 6.  However, with so many vertices and edges we cannot readily read this information from the graph, and any extension to Sudoku would be far worse.  In the next section we develop a method of organizing this information that clarifies the action of the position symmetries up to relabeling permutations.


\section{Visualizing Shidoku Symmetries with $S_4$-Nests} \label{S:vizS4}

The work in the previous section suggests that it might be helpful to consider the action of position symmetries separately from the action of relabeling symmetries.  If we group the 288 Shidoku boards into orbits under the action of $S_4$, we can restrict our attention to the action of the position symmetries on these orbits.  We will call each orbit an {\em $S_4$-nest}, and use these nests as the vertices of a simplified graph.  

We can easily define a unique representative of each $S_4$-nest by choosing the unique board in the nest whose upper-left corner is labeled in the order shown in Figure 10.

\begin{figure}[h]
$$\begin{shidoku-block}
|1 |2 |* | *|.
| 3| 4| *| *|.
| *| *| *| *|.
| *| *| *| *|.
\end{shidoku-block}$$
\caption{Representative label ordering for the upper-left block.}
\end{figure}

There are exactly twelve such representatives, which we name $A$--$L$ as shown in Figure 11.  Each of the twelve corresponding $S_4$-nests will also be denoted by $A$--$L$.

We can now form the {\em $S_4$-nest graph} whose vertices are these twelve $S_4$-nests and whose edges are given by the induced action of the position symmetries $r$, $s$, and $t$ on those $S_4$-nests.  As usual, this graph the edges represent our choice of generators for $H_4$ with the action of other group elements such as $st$ or $stst$ appearing as paths.

\begin{theorem}
Every Shidoku board is represented by exactly one of the twelve boards $A$--$L$ in Figure 11, and the group action of $H_4 =\langle r, s, t \rangle$ on this set of representatives is given as in the diagram in Figure 12.
Moreover, if $H'$ is a subgroup of $H_4=\langle r, s, t \rangle$ then $H' \times S_4$ is a complete Shidoku symmetry group if and only if the edges corresponding to a set of generators of $H'$ form a two-component graph on the vertices $A$--$L$.
\label{T:S4nests}
\end{theorem}


\begin{figure}[hb]
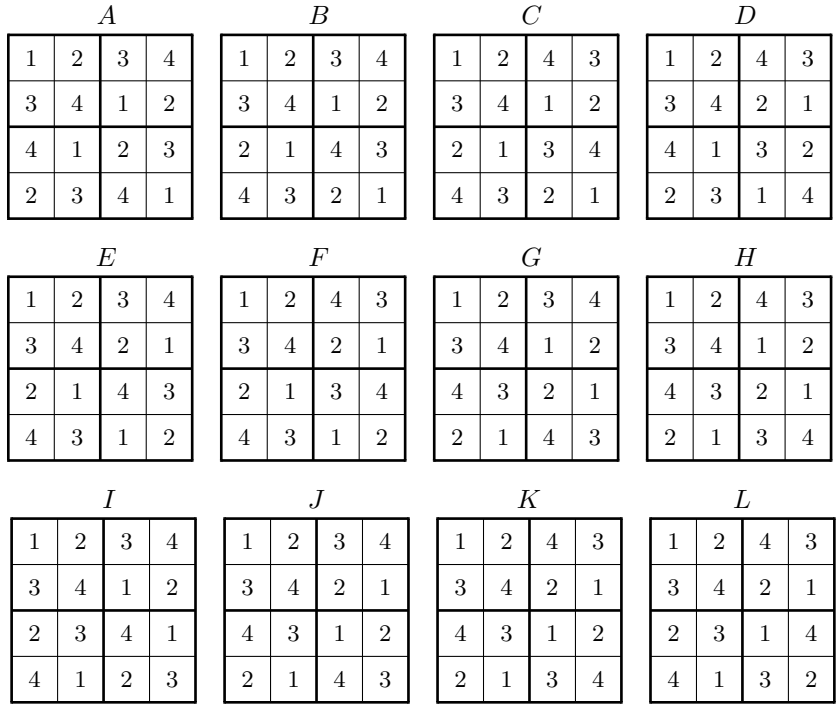

{\em \hspace{.04in} 
$A$ \hspace{.935in} $B$ \hspace{.935in} $C$ \hspace{.935in} $D$}
\vspace{.25pc}

\begin{shidoku-block}
|1 |2 |3 | 4|.
| 3| 4| 1| 2|.
| 4| 1| 2| 3|.
| 2| 3| 4| 1|.
\end{shidoku-block}
\hspace{.5pc}
\begin{shidoku-block}
|1 |2 |3 | 4|.
| 3| 4| 1| 2|.
| 2| 1| 4| 3|.
| 4| 3| 2| 1|.
\end{shidoku-block}
\hspace{.5pc}
\begin{shidoku-block}
|1 |2 |4| 3|.
| 3| 4| 1| 2|.
| 2| 1| 3| 4|.
| 4| 3| 2| 1|.
\end{shidoku-block}
\hspace{.5pc}
\begin{shidoku-block}
|1 |2 |4 | 3|.
| 3| 4| 2| 1|.
| 4| 1| 3| 2|.
| 2| 3| 1| 4|.
\end{shidoku-block}
\vspace{.5pc}

{\em \hspace{.04in} $E$ \hspace{.935in} $F$ \hspace{.935in} $G$ \hspace{.935in} $H$}
\vspace{.25pc}

\begin{shidoku-block}
|1 |2 |3 | 4|.
| 3| 4| 2| 1|.
| 2| 1| 4| 3|.
| 4| 3| 1| 2|.
\end{shidoku-block}
\hspace{.5pc}
\begin{shidoku-block}
|1 |2 |4 | 3|.
| 3| 4| 2| 1|.
| 2| 1| 3| 4|.
| 4| 3| 1| 2|.
\end{shidoku-block}
\hspace{.5pc}
\begin{shidoku-block}
|1 |2 |3 | 4|.
| 3| 4| 1| 2|.
| 4| 3| 2| 1|.
| 2| 1| 4| 3|.
\end{shidoku-block}
\hspace{.5pc}
\begin{shidoku-block}
|1 |2 |4 | 3|.
| 3| 4| 1| 2|.
| 4| 3| 2| 1|.
| 2| 1| 3| 4|.
\end{shidoku-block}
\vspace{.5pc}

{\em \hspace{.04in} $I$ \hspace{.935in} $J$ \hspace{.935in} $K$ \hspace{.935in} $L$}
\vspace{.25pc}

\begin{shidoku-block}
|1 |2 |3 | 4|.
| 3| 4| 1| 2|. 
| 2| 3| 4| 1|.
| 4| 1| 2| 3|.
\end{shidoku-block}
\hspace{.5pc}
\begin{shidoku-block}
|1 |2 |3 | 4|.
| 3| 4| 2| 1|.
| 4| 3| 1| 2|.
| 2| 1| 4| 3|.
\end{shidoku-block}
\hspace{.5pc}
\begin{shidoku-block}
|1 |2 |4 | 3|.
| 3| 4| 2| 1|.
| 4| 3| 1| 2|.
| 2| 1| 3| 4|.
\end{shidoku-block}
\hspace{.5pc}
\begin{shidoku-block}
|1 |2 |4| 3|.
| 3| 4| 2| 1|.
| 2| 3| 1| 4|.
| 4| 1| 3| 2|.
\end{shidoku-block}
\caption{The sixteen $S_4$-nest representatives.}
\end{figure}

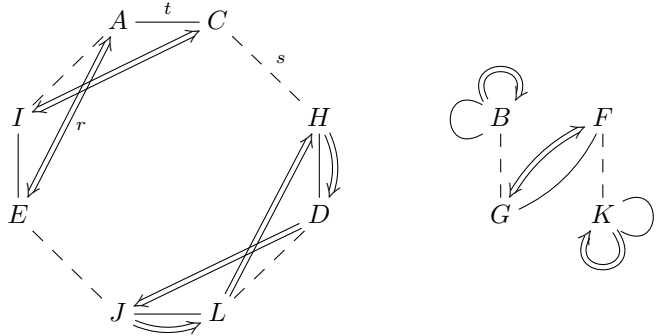
\begin{figure}[h,t]
$$\xymatrix{&A\ar@{-}[r]^t\ar@{--}[dl]\ar@{<=>}[ddl]^r&C\ar@{--}[dr]^s&&&&\\
I\ar@{-}[d]\ar@{<=>}[urr]&&&H\ar@{-}[d]\ar@/^/@{=>}[d]&&B\ar@{--}[d]\ar@(ul,dl)@{-}[]\ar@(ul,ur)@{=>}[]&F\ar@{--}[d]\ar@/^/@{-}[dl]\ar@/_/@{<=>}[dl]\\
E\ar@{--}[dr]&&&D\ar@{--}[dl]\ar@{=>}[dll]&&G&K\ar@(ur,dr)@{-}[]\ar@(dr,dl)@{=>}[]\\
&J\ar@{-}[r]\ar@/_/@{=>}[r]&L\ar@{=>}[uur]&&&&&}
$$
\caption{Group action of $H_4$ on the set of $S_4$-nests.}
\end{figure}

Notice that the $S_4$-nest graph has two components because there are two equivalence classes of Shidoku boards with respect to the Shidoku symmetry group.  Each of the $S_4$-nests contains $4!$ different Shidoku boards, and thus the larger component of this graph represents $8(4!) = 192$ Shidoku boards and the smaller component represents $4(4!) = 96$ Shidoku boards, as expected.

\begin{proof}
By construction, the action of $H_4 = \langle r, s, t \rangle$ on the nests $A$--$L$ is well-defined.  If $H'$ is some subset of $H=\langle r, s, t \rangle$ then the condition of being a complete Shidoku symmetry group is equivalent to having two orbits under the action of this group, which in turn means that the associated graph on $A$--$L$ will have two components.

Proving that the diagram in Figure 12 holds is a simple matter of applying $r$, $s$, and $t$ to each of the twelve representative boards and then, if necessary,  applying a relabeling symmetry in order to obtain one of the twelve representative boards $A$--$L$.
\end{proof}

The edge labeled $t$ between nests $A$ and $C$ indicates that there is some relabeling symmetry $\sigma \in S_4$ such that $\sigma(t(A))=C$.  Because of the way we chose our representatives $A$--$L$, the relabeling for $t$ is  always $(2\,\,\,3)$.  The position symmetry $s$, on the other hand, does not affect the first block so no relabeling is needed after applying $s$ (e.g. $s(C)=H$).  Note also that $t$ and $s$ are their own inverses and therefore their edges do not have directional arrows.

The double arrows represent application of the position rotation symmetry $r$ followed by whichever labeling is needed to get the top left square in the $1,2,3,4$ order.  This relabeling symmetry is different for different boards.  Notice that for some $S_4$-nests the action of $r$ is its own inverse, and for others it is not. 
\vspace{1pc}

\noindent
{\em Example 6: $\langle s, t \rangle$ on $S_4$-nests}

\noindent
From the $S_4$-nest graph we can immediately see that removing the edges labeled $r$ preserves the number of components, so the subgroup $\langle s, t \rangle \times S_4$ (first seen in Example 3) is a minimal complete Shidoku symmetry group, as shown in Figure 13.

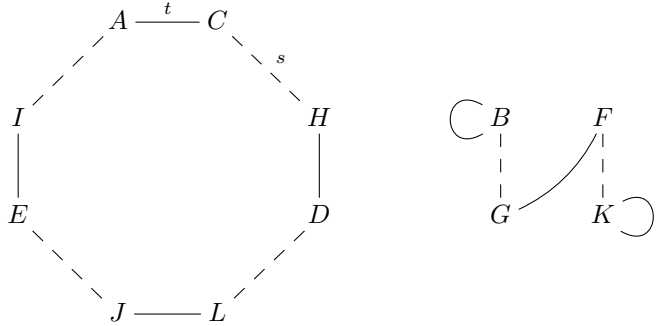
\begin{figure}[h,t]
$$\xymatrix{&A\ar@{-}[r]^t\ar@{--}[dl]&C\ar@{--}[dr]^s&&&&\\
I\ar@{-}[d]&&&H\ar@{-}[d]&&B\ar@{--}[d]\ar@(ul,dl)@{-}[]&F\ar@{--}[d]\ar@/^/@{-}[dl]\\
E\ar@{--}[dr]&&&D\ar@{--}[dl]&&G&K\ar@(ur,dr)@{-}[]\\
&J\ar@{-}[r]&L&&&&&}$$
\caption{The action of $\langle s, t \rangle$ on $S_4$-nests.}
\end{figure}

\vspace{1pc}
\noindent
{\em Example 7: $\langle r, t \rangle$ on $S_4$-nests}

\noindent
Similarly, we can see that removing the edges labeled $s$ results in a five-component graph, and thus that $\langle r, t \rangle \times S_4$ is not a complete Shidoku symmetry group, as we saw in Example 1; see Figure 14.

\begin{figure}[h,t]
$$\xymatrix{&A\ar@{-}[r]^t\ar@{<=>}[ddl]^r&C&&&&\\
I\ar@{-}[d]\ar@{<=>}[urr]&&&H\ar@{-}[d]\ar@/^/@{=>}[d]&&B\ar@(ul,dl)@{-}[]\ar@(ul,ur)@{=>}[]&F\ar@/^/@{-}[dl]\ar@/_/@{<=>}[dl]\\
E&&&D\ar@{=>}[dll]&&G&K\ar@(ur,dr)@{-}[]\ar@(dr,dl)@{=>}[]\\
&J\ar@{-}[r]\ar@/_/@{=>}[r]&L\ar@{=>}[uur]&&&&&}$$
\caption{Action of $\langle r, t \rangle$ on $S_4$-nests (too many components).}
\end{figure}
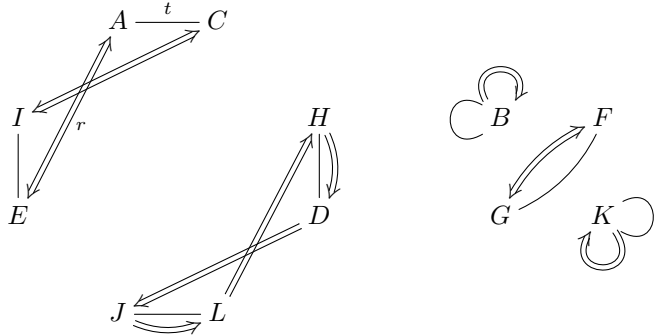
\vspace{1pc}

We can also use the $S_4$-nest graph in Figure 12 to immediately obtain the invariance data in Table 1 used for Burnside's Lemma.  For example, since the edges labeled $s$ never form loops at any $S_4$-nest, there is no Shidoku board that is invariant under the action of $s$ up to relabeling.  Similarly, there are two $S_4$-nests for which the action of $t$ forms a loop, so $2(4!)$ Shidoku boards that are invariant under the action of $t$ up to relabeling.

\section{Visualizing Shidoku Symmetries with $H_4$-Nests}

The $S_4$-nest graph can only be used to investigate Shidoku symmetry subgroups of the form $H' \times S_4$.  To investigate subgroups of the form $H_4 \times S^{\prime}$ for  subgroups $S^{\prime} \subset S_4$, we follow a dual procedure.  The first step is to partition the 288 Shidoku boards into nests according to the orbits of the action of $H_4 = \langle r, s, t \rangle$ on those boards.   We can then create a graph whose vertices correspond to nests of the group of position symmetries $H_4$ and whose edges correspond to elements of $S_4$.  The process will be similar to, but less straightforward than, the process we used to find $S_4$-nests.

Using {\em yEd} we can make a graph with 288 vertices for the Shidoku boards, and edges drawn only for position symmetries in $H_4=\langle r,s,t \rangle$ as shown in Figure 15.  From this graph we can see that there should be six $H_4$-nests, three of which are large and three of which are small.   This makes sense as some Shidoku boards are invariant under certain position symmetries (for example, some boards are preserved under transpose), so nests containing such boards should be smaller.

\begin{figure}[h,t]
\centerline{\includegraphics[height=3in, width=4.3in]{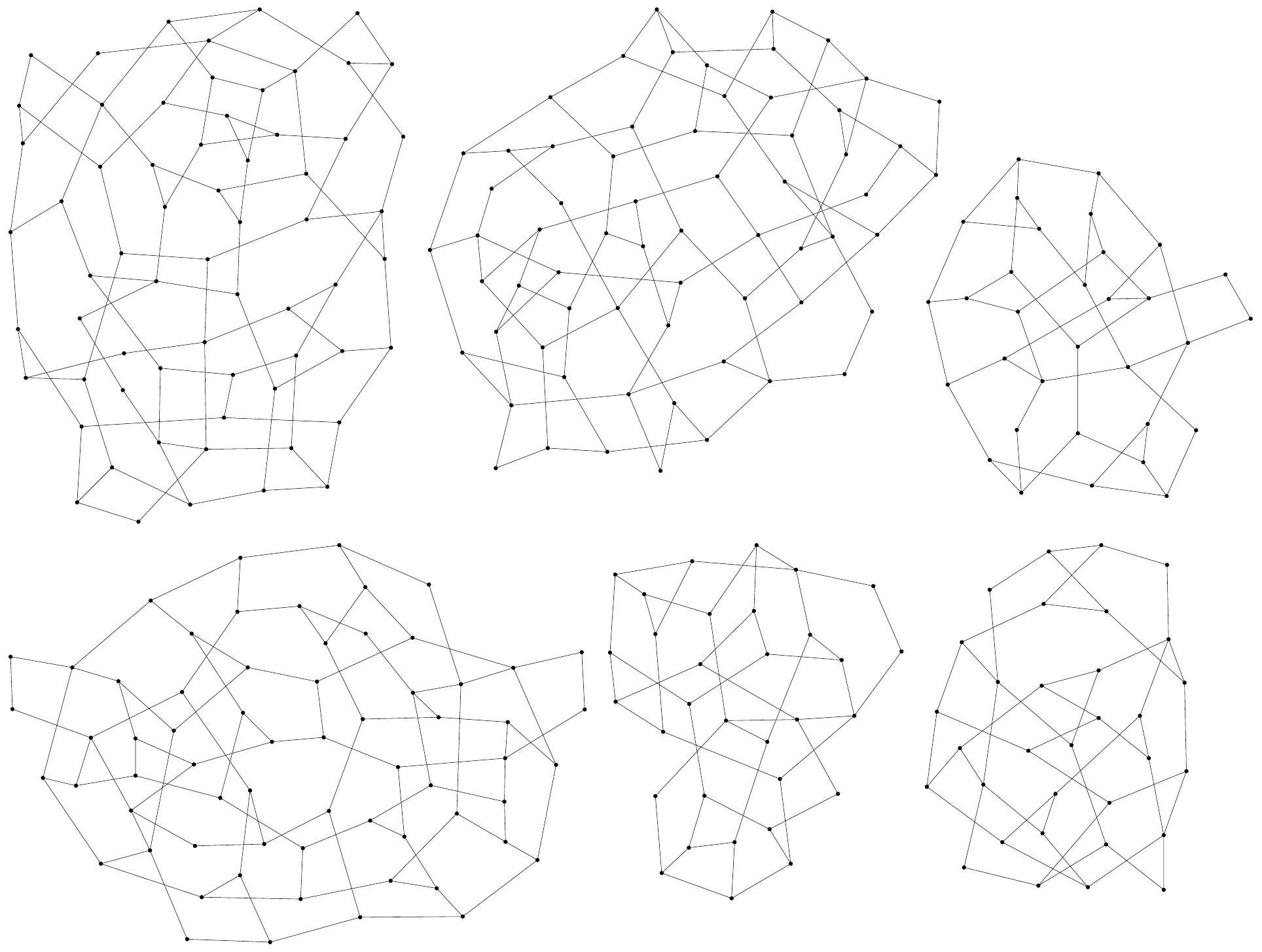}}
\caption{$H_4$-nests arising from the action of $H_4$ on $\mathcal{S}$.}
\end{figure}

We must now find a way to choose a unique representative Shidoku board for each $H_4$-nest.  When we partitioned the 288 Shidoku boards into $S_4$-nests we chose one unique representative $A$--$L$ for each $S_4$-nest by fixing the form of the upper-left block of the board.  This very neatly used all of the $S_4$ part of the action of $H_4 \times S_4$.  With the $S_4$-nests it is very easy to take any Shidoku board and immediately find its unique representative, and it is also easy to show that there are twelve such representatives and thus twelve $S_4$-nests, each containing $4!$ Shidoku boards.

Similarly, the process of finding such a representative for the $H_4$-nests should use all of the $H_4$ part of the action of $H_4 \times S_4$, and it should be easy to take any Shidoku board and find its unique representative.  The difficulty lies in the fact that we must find this representative by position symmetries only.  

\begin{theorem}
Each $H_4$-nest has a unique representative of the form shown in Figure 16, 
where $a\leq b$ and $c<d$.  Moreover, there are six $H_4$-nests, three of size 32 and three of size 64.
\label{T:Hrepresentatives}
\end{theorem}

\begin{figure}[h,t]
$$\begin{shidoku-block}
|1 |c |*| *|.
| d| a| 1| *|.
| *| 1| b| *|.
| *| *| *| 1|.
\end{shidoku-block}$$
\caption{Form of $H_4$-nest representatives.}
\end{figure}

\begin{proof}
Given any Shidoku board we can perform a sequence of row and column swaps to obtain a board whose configuration of 1s is the same as in Figure 16.  A calculation similar to the one commonly used for enumerating Shidoku boards (\cite{Arnold}, \cite{Taalman}) shows that there are $18$ Shidoku boards that have this particular configuration of 1s.

Now via the possible application of the $180^{\circ}$ rotation symmetry $r^2$ we can force the inequality $a \leq b$ in the 6th and 11th cells of the board as ordered in Figure 8.   Note that this rotation does not change the configuration of 1s we had already established.  The equality $a=b$ holds for Type 1 boards (see Figure 2).  In each of the three cases $a=b=2$, $3$, and $4$ there are two possible Shidoku completions.  We have strict inequality $a<b$ for Type 2 boards, and it is a simple calculation to verify that in each of the three possible cases $2<3$, $2<4$, and $3<4$ we will also have two Shidoku completions.  Therefore there are six Shidoku boards that satisfy both the 1s configuration and the $a = b$ condition, and six Shidoku boards that satisfy both the 1s configuration and the $a<b$ condition.  

Finally, by applying the transpose symmetry $t$ if necessary, we can require that $c<d$ in the 2nd and 5th cells.  It is a simple matter to check that any board with configuration of 1s as given in Figure 16 that satisfies $a \leq b$, $c<d$ in the relevant cells can have only one possible completion.   Since there are six ways to select $a$ and $b$ from the set $\{2,3,4\}$ with $a \leq b$, and for each of these choices there is only one way to complete the first block with $c<d$, there are six $H_4$-nest representatives.  

The six $H_4$-nest representatives constructed above are shown in Figure 17 and labeled $a$--$f$ based on their lexicographical order.  

There are two possible sizes of nests.   Nests $b$, $d$, and $e$ contain $64$ Shidoku boards.  The remaining nests $a$, $c$, and $f$ are rotationally symmetric and therefore contain only $32$ Shidoku boards. 
\end{proof}


\begin{figure}[h]
{\em \hspace{-.1in} $a$ \hspace{1.23in} $b$ \hspace{1.23in} $c$}
\vspace{.5pc}

\hfill
\begin{shidoku-block}
|1 |2 |3 | 4|.
| 3| 4| 1| 2|. 
| 2| 1| 4| 3|.
| 4| 3| 2| 1|.
\end{shidoku-block}
\hspace{2pc}
\begin{shidoku-block}
|1 |2 |3 | 4|.
| 4| 3| 1| 2|.
| 2| 1| 4| 3|.
| 3| 4| 2| 1|.
\end{shidoku-block}
\hspace{2pc}
\begin{shidoku-block}
|1 |2 |4 | 3|.
| 4| 3| 1| 2|.
| 2| 1| 3| 4|.
| 3| 4| 2| 1|.
\end{shidoku-block}
\hfill
\phantom{x}
\vspace{.5pc}

{\em \hspace{.01in} $d$ \hspace{1.23in} $e$ \hspace{1.23in} $f$}
\vspace{.5pc}

\begin{shidoku-block}
|1 |3 |2| 4|.
| 4| 2| 1| 3|.
| 3| 1| 4| 2|.
| 2| 4| 3| 1|.
\end{shidoku-block}
\hspace{2pc}
\begin{shidoku-block}
|1 |3 |4 | 2|.
| 4| 2| 1| 3|. 
| 2| 1| 3| 4|.
| 3| 4| 2| 1|.
\end{shidoku-block}
\hspace{2pc}
\begin{shidoku-block}
|1 |3 |4 | 2|.
| 4| 2| 1| 3|.
| 3| 1| 2| 4|.
| 2| 4| 3| 1|.
\end{shidoku-block}
\hfill
\caption{The six $H_4$-nest representatives.}
\end{figure}

Note that we never used $r$ in the proof of Theorem~\ref{T:Hrepresentatives}; we used only $r^2$ and row and column swaps.   In fact, the position symmetry $r$ does not increase the size of the orbits and $\langle r^2, s, t \rangle$-nests are the same as the $\langle r, s, t \rangle=H_4$-nests.  This is despite the fact that the two groups  are not equal (for example, $\langle r^2,s,t\rangle$ does not contain the position symmetries of band or pillar swaps).  This will be important later in the section.

We can now form the {\em $H_4$-nest graph} whose vertices are the six $H_4$-nests $a$--$f$ and whose edges are given by the action of the relabeling symmetries $(12)$, $(23)$, $(34)$, and $(14)$.  Once again we use the same names for the $H_4$-nests as for their representatives and choose an overdetermined set of generators for $S_4$ for symmetry.  

\begin{theorem}
Every Shidoku board is represented by one of the six boards $a$--$f$ in Figure 17, and the group action of $S_4=\langle (12), (23), (34), (14) \rangle$ on this set of representatives is given as in the diagram in Figure 18.
Moreover, if $S'$ is a subgroup of $S_4$ and $H_4=\langle r, s, t \rangle$, then $H_4 \times S'$ is a complete Shidoku symmetry group if and only if the edges corresponding to a set of generators from $S'$ form a two-component graph on the vertices $a$--$f$.
\label{T:Hnests}
\end{theorem}

\begin{figure}[h,t]
$$\xymatrix{*++[o][F-]{a}\ar@(ul,dl)@{-}[]_{(14)}\ar@(ul,ur)@{-}[]^{(23)}\ar@/^/@{-}[rr]^{(34)}\ar@/_/@{-}[rr]|{(12)}&&*++[o][F-]{c}\ar@/^/@{-}[dl]^{(14)}\ar@/_/@{-}[dl]|{(23)}&&&*++[o][F-]{b}\ar@(ul,dl)@{-}[]_{(14)}\ar@(ul,ur)@{-}[]^{(23)}\ar@/^/@{-}[rr]^{(34)}\ar@/_/@{-}[rr]|{(12)}&&*++[o][F-]{d}\ar@/^/@{-}[dl]^{(14)}\ar@/_/@{-}[dl]|{(23)}\\
&*++[o][F-]{f}\ar@(ul,dl)@{-}[]_{(12)}\ar@(dl,dr)@{-}[]_{(34)}&&&&&*++[o][F-]{e}\ar@(ul,dl)@{-}[]_{(12)}\ar@(dl,dr)@{-}[]_{(34)}}$$
\caption{Group action of $S_4$ on the set of $H_4$-nests.}
\end{figure}
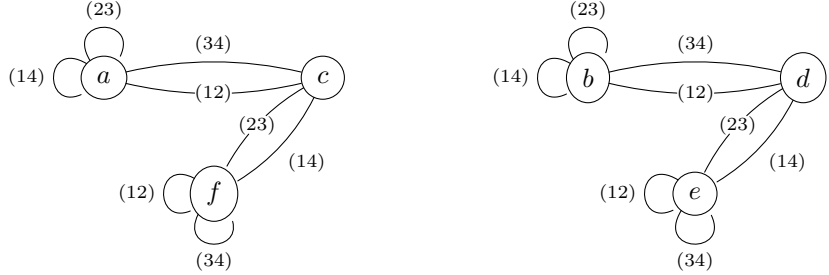

The proof of Theorem~\ref{T:Hnests} is entirely similar to the proof to Theorem~\ref{T:S4nests}, with the details of the $H_4$-nest representatives established by Theorem~\ref{T:Hrepresentatives}.  We give two examples:  applying the transposition $(34)$ to the representative board $a$ gives us the representative board $c$ directly; so $(34)(a) = c$.   Applying $(23)$ to the board $a$ gives us us a board which after transposing is once again equal to $a$ or $t((23)(a))=a$.

Just as the $S_4$-nest graph allowed us to immediately identify ways to reduce the size of $H_4 \times S_4$ by eliminating position symmetries in $H_4$, this dual $H_4$-nest graph allows us to immediately identify redundant relabeling symmetries.  
\vspace{1pc}

\noindent
{\em Example 8: $\langle (12), (23) \rangle$ on $H_4$-nests}

\noindent
Removing the transpositions $(14)$ and $(34)$ from the $H_4$-nest graph in Figure 18 preserves the number of components and therefore the 384-element subgroup $H_4 \times \langle (12),(23) \rangle$ is a full Shidoku symmetry group, although not minimal; see Figure 19.

\begin{figure}[h,t]
$$\xymatrix{*++[o][F-]{a}\ar@(ul,ur)@{-}[]^{(23)}\ar@/_/@{-}[rr]|{(12)}&&*++[o][F-]{c}\ar@/_/@{-}[dl]|{(23)}&&&*++[o][F-]{b}\ar@(ul,ur)@{-}[]^{(23)}\ar@/_/@{-}[rr]|{(12)}&&*++[o][F-]{d}\ar@/_/@{-}[dl]|{(23)}\\
&*++[o][F-]{f}\ar@(ul,dl)@{-}[]_{(12)}&&&&&*++[o][F-]{e}\ar@(ul,dl)@{-}[]_{(12)}\ar@(dl,dr)@{-}[]_{(34)}}$$
\caption{The action of $\langle (12),(23) \rangle$ on $H_4$-nests.}
\end{figure}
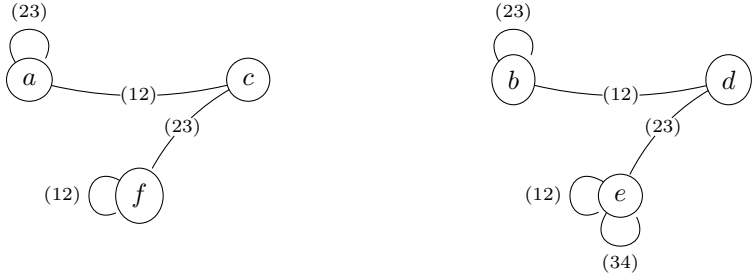

\noindent {\em Example 9: $\langle (123) \rangle$ on $H_4$-nests}

\noindent There is a smaller class of subgroups of $S_4$ which also produce complete Shidoku symmetry groups: the ones of order three.  The example of $\langle (123)\rangle$ is given in Figure 20.

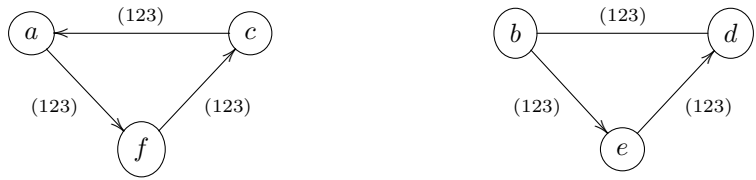
\begin{figure}[h,t]
$$\xymatrix{*++[o][F-]{a}\ar[dr]_{(123)}&&*++[o][F-]{c}\ar[ll]_{(123)}&&&*++[o][F-]{b}\ar[dr]_{(123)}&&*++[o][F-]{d}\ar@{-}[ll]_{(123)}\\
&*++[o][F-]{f}\ar[ur]_{(123)}&&&&&*++[o][F-]{e}\ar[ur]_{(123)}&
}$$
\caption{Action of $\langle (123) \rangle$ on $H_4$-nests.}
\end{figure}

Note that in this case our choice of generator $(123)$ forces us to use a different edge set than we used in the $H_4$-nest graph shown in Theorem~\ref{T:Hnests}. 
With the vertices as $H_4$-nests, this does not produce a minimal Shidoku symmetry group as $|H_4\times \langle(123)\rangle|=384$.  However, as noted in Remark 4, $H_4$-nests are the same as $\langle r^2,s,t\rangle$-nests, and the order of $\langle r^2,s,t\rangle\times\langle(123)\rangle$ is $192$ as desired.  
This is the same Shidoku symmetry group we found in Example 4 using {\em GAP} and {\em yEd}, but it is not the only example of a complete minimal Shidoku symmetry group with $\langle(123)\rangle$ for its relabeling symmetries.  In Example 2 we saw that $\langle r,s\rangle\times\langle(123)\rangle$ form a complete minimal Shidoku symmetry group.  In fact, $H_4$-nests are {\em not} the same as $\langle r,s\rangle$-nests even though the graph for $\langle r,s\rangle$-nests under the action of $(123)$ looks identical that in Figure 20.  The $\langle r,s\rangle$-nests simply contain different boards than the $H_4$-nests.

\section{Future Directions for Sudoku Symmetry Groups}

Using $H_4$-nests or $S_4$-nests we can quickly see by hand what previously required computing with {\em yEd} and {\em GAP} software.

In \cite{Russell}, Russell and Jarvis use a conjugacy class table to find the number of equivalence classes of $9 \times 9$ Sudoku boards using Burnside's Lemma.  {\em GAP} quite easily computes the full Sudoku symmetry group generated by all of the position symmetries described in \cite{Felgenhauer} and \cite{Russell}.  The size of the full position symmetry group $H_9$ is 3,359,232. The full Sudoku symmetry group $G_9 = H_9 \times S_9$ has order 1,218,998,108,160.  Since this is the full Sudoku symmetry group, it is complete.  But is it minimal?  

Jarvis and Russell compute that there are 5,472,730,538 equivalence classes of Sudoku boards. They use a computer to find the number of invariant boards for each conjugacy class.  The number of conjugacy classes of the full Sudoku symmetry group is intractable to work with by hand, but with nests and reduced complete Sudoku symmetry groups there may be a less computational way to compute this number of equivalence classes.  However, there are important differences between Shidoku and Sudoku.  In the $9 \times 9$ Sudoku case, the size of the full symmetry group is smaller than the number of Sudoku boards, while in the $4 \times 4$ Shidoku case, the size of the full symmetry group is significantly larger than the size of the largest orbit.  For Sudoku we do not even know the size of the largest orbit.  It is even possible that the full symmetry group for Sudoku is actually minimal. 

Some problems lie between Shidoku and Sudoku, for example various modular or magic Sudoku-type boards, the $6 \times 6$ Roku-Doku boards, and boards whose nonlinear regions are built from polyominoes (for example the CrossDoku variation from \cite{gerechte}).  Applying orbit graphs and nest techniques could lead to finding minimal complete symmetry groups for these intermediate examples, and bring insight and understanding to the larger Sudoku symmetry group.



\begin{thebibliography}{99}

\bibitem{Arnold} E. Arnold, S. Lucas \& L. Taalman, Gr\"obner basis representations of Sudoku, \textit{College Mathematics Journal}, \textbf{41}(2), 2010, 101--112.

\bibitem{gerechte}R.A.Bailey, P.Cameron \& R. Connelly, Sudoku, gerechte designs, resolutions, affine space, spreads, reguli, and Hamming codes, \textit{Sudoku American Mathematical Monthly}, \textbf{115}(5), 2008, 383--404(22).

\bibitem{GAP} The GAP Group, GAP -- Groups, Algorithms, and Programming, Version 4.4.12, 2008. {\tt (http://www.gap-system.org)}

\bibitem{Felgenhauer} B. Felgenhauer \& F. Jarvis, Mathematics of Sudoku I, \textit{Mathematical Spectrum}, \textbf{39}, 2006, 15--22.

\bibitem{lorch} C. Lorch \& J Lorch, Enumerating small Sudoku puzzles in a first abstract algebra course, {\em PRIMUS} {\bf 18}, 2008, 149--158.

\bibitem{Newton} P. Newton \& S. DeSalvo, The Shannon entropy of Sudoku matrices, \textit{Proc. Royal Soc. A}, \textbf{466}, 2010, 1957--1975. 

\bibitem{Matlab} MATLAB, Natick, Massachusetts, The MathWorks Inc., 2003.

\bibitem{nofrills} P. Riley \& L. Taalman, {\em No-Frills Sudoku}, Puzzle Wright Press, Sterling, 2011.

\bibitem{Russell} E. Russell \& F. Jarvis, Mathematics of Sudoku II, \textit{Mathematical Spectrum}, \textbf{39}, 2006, 54--58.

\bibitem{Sato} Y. Sato, S. Inoue, A. Suzuki, K. Nabeshima, \& K. Sakai, Boolean Gr\"obner bases, \textit{Journal of Symbolic Computation}, \textbf{46}(5), 2011, 622--632.

\bibitem{Taalman} L. Taalman, Taking Sudoku seriously, \textit{Math Horizons}, September 2007, 5--9.

\bibitem{Wilson} R. Wilson, The Sudoku Epidemic, \textit{FOCUS}, January 2006, 5--7.

\bibitem{yed} yEd version 3.7.02, T\"ubingen, Gemany: yWorks GmbH, 2011.

\end{thebibliography}
\end{document}